\documentclass[11pt]{amsart}
\usepackage{fullpage}
\usepackage {amsmath}
\usepackage {amsthm}
\usepackage {enumerate}
\usepackage{amssymb}
\usepackage{amsfonts}
\usepackage{anysize}
\usepackage{setspace}

\marginsize{3.15cm}{3.15cm}{1.6cm}{3cm}

\newcommand{\N}{\mathbb N}

\newcommand{\R}{\mathbb R}

\numberwithin{equation}{section}

\theoremstyle{plain}

\newtheorem{theorem}[equation]{Theorem}
\newtheorem{proposition}[equation]{Proposition}
\newtheorem{corollary}[equation]{Corollary}
\newtheorem{lemma}[equation]{Lemma}
\newtheorem{question}[equation]{Question}
\newtheorem{problem}[equation]{Problem}

\theoremstyle{definition}

\theoremstyle{definition}

\theoremstyle{definition}

\theoremstyle{definition}
\newtheorem{example}[equation]{Example}

\theoremstyle{definition}

\author{{\bfseries S. Garc\'ia-Ferreira}}
\address{Centro de Ciencias Matem\'aticas\\
         Universidad Nacional Aut\'onoma de M\'exico\\
				 Campus Morelia\\
         Apartado Postal 61-3, Santa Mar\'ia, 58089, Morelia, Michoac\'an, M\'exico.}
\email{sgarcia@matmor.unam.mx}

\title{\scshape\bfseries Connectedness like properties on the hyperspace of convergent sequences}

\author{{\bfseries R. Rojas-Hern\'andez}}
\email{satzchen@yahoo.com.mx}



\date{}

\thanks{Research of the first-named author was supported
by  CONACYT grant no. 176202 and PAPIIT grant no. IN-101911.}

\begin{document}

\begin{abstract}
This paper is a continuation of the work done in \cite{sal-yas} and \cite{may-pat-rob}. We deal with the Vietoris hyperspace of all nontrivial convergent sequences $\mathcal{S}_c(X)$ of a space $X$. We answer some questions in \cite{sal-yas} and generalize several results in \cite{may-pat-rob}. We prove that: The connectedness of $X$ implies the connectedness of $\mathcal{S}_c(X)$; the local connectedness of $X$ is equivalent to the local connectedness of $\mathcal{S}_c(X)$; and the path-wise connectedness of $\mathcal{S}_c(X)$ implies the path-wise connectedness of $X$. We also show that the space of nontrivial convergent sequences on the Warsaw circle has $\mathfrak{c}$-many path-wise connected components, and provide a dendroid with the same property.
\end{abstract}

\maketitle

\section{Introduction}

 The letters $\mathbb{R}$,  $\N$  and $\mathbb{Z}$ will denote the real numbers,   natural numbers and the integers, respectively. In particular, $\N^+ = \N \setminus\{0\}$. The Greek letter $\omega$ stands for the first infinite cardinal number and $\mathfrak{c}$ will denote the cardinality of the continuum.

\medskip

For a topological space $X$, $\mathcal{CL}(X)$ will denote the set of all nonempty closed subsets of $X$. For a nonempty family $\mathcal{U}$ of subsets of $X$ let
\begin{center}
$\left\langle \mathcal{U} \right\rangle = \{F \in \mathcal{CL}(X) : F \subseteq \bigcup{\mathcal{U}} \textnormal{ and } F \cap U \not= \emptyset\}.$
\end{center}
If $\mathcal{U} = \{U_1,\ldots,U_n\}$, in some cases,  $\left\langle \mathcal{U} \right\rangle$ will be denoted by  $\left\langle U_1,\ldots,U_n \right\rangle$. A base for the \textit{Vietoris topology} on $\mathcal{CL}(X)$ is the family of all sets of the form $\left\langle \mathcal{U}\right\rangle$, where $\mathcal{U}$ runs over all nonempty finite families of nonempty open subsets of $X$. In the sequel, any subset $\mathcal{D} \subseteq \mathcal{CL}(X)$ will carry the relative Vietoris topology as a subspace of $\mathcal{CL}(X)$. Given $\mathcal{D} \subseteq \mathcal{CL}(X)$ and a nonempty family $\mathcal{U}$ of subsets of $X$ we let $\left\langle \mathcal{U}\right\rangle_\mathcal{D} = \left\langle \mathcal{U}\right\rangle \cap \mathcal{D}$. For simplicity, if there is no possibility of confusion we simply write  $\left\langle \mathcal{U}\right\rangle$ instead of $\left\langle \mathcal{U}\right\rangle_\mathcal{D}$.
All topological notions whose definition is not included here should be understood as in \cite{eng}.

\medskip

Some of the most studied hyperspaces on a space $X$ have been
$$
\mathcal{K}(X) = \{K \in \mathcal{CL}(X) : K \textnormal{ is compact}\} \textnormal{ and } \mathcal{F}(X) = \{F \in \mathcal{CL}(X) : F \textnormal{ is finite}\},
$$
see for instance the survey paper \cite{gu}.
In the article \cite{jvpp}, the authors considered the hyperspace consisting of all finite subsets together with all  the Cauchy sequences  without limit point of a metric space.
In a different context,  the authors of the paper \cite{peng} consider the set $\mathcal{F}_S(X)$ of all convergent sequences of a space $X$, and study the existence of a metric $d$ on the set $X$ such that $d$ metrizes all subspaces of $X$ which belong to $\mathcal{F}_S(X)$, that is, the restriction of $d$ to $A$ generates the subspace topology on $A$ for every $A \in \mathcal{F}_S(X)$ (these kind of problems have been analyzed in \cite{arh}).

\medskip

 Another hyperspace that has been recently introduced in \cite{sal-yas} and studied in \cite{may-pat-rob} is the hyperspace of nontrivial convergent sequences which is defined as follows:

\smallskip

 Given a space $X$, a \textit{nontrivial convergent sequence} of $X$ is a subset $S \subseteq X$ such that  $\left|S\right| = \omega$, $S$ has a unique non isolated point $x_S$ and $\left|S \setminus U\right| < \omega$ for each neighborhood $U$ of $x_S$. With this notion we define
$$\mathcal{S}_c(X) = \{S \in \mathcal{CL}(X) : S \textnormal{ is a nontrivial convergent sequence}\}.$$
It was pointed out in \cite{may-pat-rob} that the family of all subsets of $\mathcal{S}_c(X)$ of the form $\left\langle \mathcal{U}\right\rangle$, where $\mathcal{U}$ is a finite family of pairwise disjoint subsets of $X$,  is a base for the Vietoris topology on $\mathcal{S}_c(X)$. We will refer to  this family as the {\it canonical basis} of $\mathcal{S}_c(X)$ and  its elements will be name {\it canonical open sets}.

\medskip

Along this paper we will deal with hyperspaces of non trivial convergent sequences. By this reason, it is natural to consider only  Fr\'echet-Urysohn spaces that are infinite and non-discrete.

\medskip

A fundamental task in the study of the hyperspace $\mathcal{S}_c(X)$ is to determine its topological relationship with  the base space $X$ and vise versa.
The connection  between the connectedness in $X$ and $\mathcal{S}_c(X)$  was studied in \cite{sal-yas}, where it was proved that $\mathcal{S}_c(X)$ is connected provided that $X$ is a path-wise connected space and the authors formulate the question (\cite[Q. 2.16]{sal-yas}) whether or not  the connectedness of $X$ implies the connectedness of $\mathcal{S}_c(X)$. Very recently, it was proved in \cite{may-pat-rob}  that if a connected space $X$ has finitely many path components or has a dense path component, then $\mathcal{S}_c(X)$ is connected. It is important to mention that, as it was proved in \cite{may-pat-rob}, the connectedness of $X$ also implies the connectedness of $\mathcal{S}_c(X)$. In contrast with the case of the connectedness, an example of a path-wise connected space such that $\mathcal{S}_c(x)$ is not path-wise connected was provided in \cite{sal-yas}. In fact, the authors showed that  the hyperspace non trivial convergent sequences over the Warsaw circle is not path-wise connected.

\medskip

 We shall show that $\mathcal{S}_c(X)$ is connected whenever $X$ is connected, this answers Question 2.16 from \cite{sal-yas}. Hence, we conclude that  $\mathcal{S}_c(X)$ connected iff $X$ connected. To end the second section we prove that $\mathcal{S}_c(X)$ is locally connected iff $X$ is locally connected.
The path-wise connectedness is considered in the third section addressing to  Question 2.14 from \cite{sal-yas} which asked
  if $\mathcal{S}_c(X)$ is path-wise connectedness, must $X$ be  path-wise connectedness? A remarkable related result was obtained in \cite{may-pat-rob}: for every metric space $X$, if $\mathcal{S}_c(X)$ nonempty and path-wise connected, then so is $X$. We prove that, in general, the path-wise connectedness of $\mathcal{S}_c(X)$ implies the path-wise connectedness of $X$.
  Answering a question from \cite[Q. 2.9]{sal-yas}, an example of a dendroid such that its hyperspace non trivial convergent sequences  is not path-wise connected was given in  \cite{may-pat-rob}. We prove that this example and the one given in \cite[Ex. 2.8]{sal-yas} both have $\mathfrak{c}$-many path-wise connected components.

\section{Connectedness and local connectedness in $\mathcal{S}_c(X)$}

To begin this section we shall prove that  $\mathcal{S}_c(X)$ is connected iff $X$ is connected. But first, we establish a basic lemma.

\begin{lemma}\label{LUA}
If $\left\langle \mathcal{U}\right\rangle$ is a nonempty canonical open set of $\mathcal{S}_c(X)$, then $\bigcup \left\langle \mathcal{U}\right\rangle = \bigcup \mathcal{U}$. Moreover, if $O$ is a nonempty open subset in $\mathcal{S}_c(X)$, then $\bigcup O$ is open in $X$.
\end{lemma}

\proof
Let $\left\langle \mathcal{U}\right\rangle$ be a nonempty canonical open set of $\mathcal{S}_c(X)$. Fix $S \in \left\langle \mathcal{U}\right\rangle$. Observe that  $x \in S_x = S \cup \{x\} \in \left\langle \mathcal{U}\right\rangle$ for all $x \in \bigcup\mathcal{U}$. It then follows that $\bigcup \mathcal{U} \subseteq \bigcup\{S_x : x \in \bigcup\mathcal{U}\} \subseteq \bigcup \left\langle \mathcal{U}\right\rangle \subseteq \bigcup \mathcal{U}$. So, $\bigcup \left\langle \mathcal{U}\right\rangle = \bigcup \mathcal{U}$.

Now assume that $O$ is a nonempty open subset of $\mathcal{S}_c(X)$. Pick $x \in \bigcup O$ and choose $S \in O$ so that $x \in S$ and a canonical open set $\left\langle \mathcal{V}\right\rangle$ of $\mathcal{S}_c(X)$ satisfying $S \in \left\langle \mathcal{V}\right\rangle \subseteq O$. Then, we have that  $x \in S \subseteq \bigcup\left\langle \mathcal{V}\right\rangle \subseteq \bigcup O$. By the first part, we conclude that $\bigcup\left\langle \mathcal{V}\right\rangle$ is open. Therefore, the set $\bigcup O$ is open.
\endproof

\begin{theorem}\label{SCC}
If $X$ is connected, then $\mathcal{S}_c(X)$ is connected.
\end{theorem}

\proof
For each $S \in \mathcal{S}_c(X)$ and $x \in X$, we define
\begin{itemize}
\item $\mathcal{G}(S) = \{S \cup F : F \in \mathcal{F}(X)\}$,

\item $\mathcal{H}(S) = \{T \in \mathcal{S}_c(X) : S \subseteq T \}$,

\item $\mathcal{M}(S) = \{T \in \mathcal{S}_c(X) : x_S = x_T\}$, and

\item $\mathcal{N}(x) = \{S \in \mathcal{S}_c(X) : x \in S\}$.\medskip
\end{itemize}
We list some basic properties of these sets:
\begin{enumerate}
\item If  $T \in \mathcal{M}(S)$, then $x_S = x_T$  and $\mathcal{H}(T) \subseteq  \mathcal{M}(S)$,

\item $\mathcal{M}(S)$ is dense in $\mathcal{N}(x_S)$,  and

\item  $\mathcal{G}(S) \subseteq \mathcal{H}(S) \subseteq \textnormal{cl}(\mathcal{G}(S))$.
 \end{enumerate}

 \medskip

\textbf{Claim.} If $S \in \mathcal{S}_c(X)$ and $x \in X$, then the sets $\mathcal{G}(S)$, $\mathcal{H}(S)$, $\mathcal{M}(S)$ and $\mathcal{N}(x)$ are connected.\medskip

    \textit{Proof of the Claim}. Since $X$ is connected we know, by Theorem 4.10 of \cite{emi}, that the hyperspace $\mathcal{F}(X)$ is connected. Let  $e : \mathcal{F}(X) \to \{S\} \times \mathcal{F}(X)$ be the map given by $g(F) = (S,F)$ and let $u : \{S\} \times \mathcal{F}(X) \to \mathcal{S}_c(X)$ be the map defined by $u(S,F) = S \cup F$. It is not hard to verify that $e$ is an embedding and $u$ is a continuous map. Hence, the composition $u \circ e$ is continuous. In addition, $u \circ e (\mathcal{F}(X)) = \mathcal{G}(S)$. Thus the connectedness of $\mathcal{G}(S)$ is obtained via the map $u \circ e$ from the connectedness of $\mathcal{F}(X)$. Hence and by clause $(3)$ we obtain that $\mathcal{H}(S)$ is connected. Next we shall verify that $\mathcal{M}(S)$ is connected.

    \medskip

    Assume on the contrary that there exist two nonempty  subsets $O_1$ and $O_2$ of $\mathcal{S}_c(X)$ which disconnect $\mathcal{M}(S)$. For each $i \in \{1,2\}$ choose $S_i \in O_i$ and a canonical open set $\left\langle \mathcal{U}_i\right\rangle$ of $\mathcal{S}_c(X)$ such that $S_i \in \left\langle \mathcal{U}_i\right\rangle \subseteq O_i$. Since $\mathcal{H}(S_i)$ is connected, we must have $\mathcal{H}(S_i) \subseteq O_i$. It follows that $S_1 \cup S_2 \subseteq O_1 \cap O_2$, a contradiction. Hence $\mathcal{M}(S)$ is connected. Finally we will prove that if $x \in X$, then $\mathcal{N}(x)$ is connected. Choose $S \in \mathcal{S}_c(X)$ converging to $x$. By clause $(2)$ and the connectedness of $\mathcal{M}(S)$,  we can conclude that  $\mathcal{N}(x)$ is connected. This ends the proof of the Claim.\medskip

Now we will prove that the space $\mathcal{S}_c(X)$ is connected. Assume that  $\mathcal{S}_c(X)$ is not connected. Then we can find two nonempty open sets $O_1$ and $O_2$ in $\mathcal{S}_c(X)$ which disconnect $\mathcal{S}_c(X)$. For every $S \in O_i$ choose a canonical open set $\left\langle \mathcal{U}_S\right\rangle$ of $\mathcal{S}_c(X)$ such that $S \in \left\langle \mathcal{U}_S\right\rangle \subseteq O_i$, for $i = 1,2$. Let $W_i = \bigcup\{\bigcup \mathcal{U}_S: S \in O_i\}$, for $i = 1,2$. According to Lemma \ref{LUA} the sets $W_1$ and $W_2$ are open. In addition, $X = W_1 \cup W_2$ because of $X$ is crowded. Since the space $X$ is connected, we must have $W_1 \cap W_2 \not= \emptyset$.  Fix $x \in W_1 \cap W_2$, then $x \in \bigcup\mathcal{U}_{S_1} \cap \bigcup\mathcal{U}_{S_2}$ for some $S_1 \in O_1$ and $S_2 \in O_2$. Choose $T_1 = S_1 \cup\{x\} \in \left\langle \mathcal{U}_{S_1}\right\rangle \subseteq O_1$ and $T_2 =  S_2 \cup\{x\} \in \left\langle \mathcal{U}_{S_2}\right\rangle \subseteq O_2$. Then $T_1,T_2 \in \mathcal{N}(x)$ and  $\mathcal{N}(x) \subseteq \mathcal{S}_c(X) \subseteq O_1 \cup O_2$, the open sets $O_1$ and $O_2$ disconnect $\mathcal{N}(x)$, which is a contradiction. Therefore,  $\mathcal{S}_c(X)$ is connected.
\endproof

It was proved in \cite[Theorem 5.10]{may-pat-rob} that $X$ must be connected whenever $\mathcal{S}_c(X)$ is connected. According to that result and Theorem  \ref{SCC}, we get the following corollary.

\begin{corollary}\label{SCCc}
The space $\mathcal{S}_c(X)$ connected iff $X$ connected.
\end{corollary}

\begin{proposition}
If $X$ is locally connected, then $\mathcal{S}_c(X)$ is locally connected.
\end{proposition}

\proof
Pick $S \in \mathcal{S}_c(X)$ and a neighborhood $O$ of $S$ in $\mathcal{S}_c(X)$. Choose a canonical open set $\left\langle \mathcal{U}\right\rangle$ of $\mathcal{S}_c(X)$ such that $S \in  \left\langle \mathcal{U}\right\rangle \subseteq O$. Since $X$ is locally connected we can find a canonical open set $\left\langle \mathcal{V}\right\rangle$ of $\mathcal{S}_c(X)$ such that $S \in  \left\langle \mathcal{V}\right\rangle \subseteq \left\langle \mathcal{U}\right\rangle$ and each member of $\mathcal{V}$ is connected. We will prove that $\left\langle \mathcal{V}\right\rangle_{\mathcal{F}(X)}$ is connected. Since $V$ is connected we know that $\mathcal{F}(V)$ is connected, for each $V \in \mathcal{V}$. Assume that $\mathcal{V} = \{V_1,\ldots,V_n\}$. Then $\prod\{\mathcal{F}(\mathcal{V}_i) : i = 1,\ldots,n\}$ is connected. Let $u : \prod\{\mathcal{F}(\mathcal{V}_i) : i = 1,\ldots,n\} \to \mathcal{F}(X)$ the map given by $u(F_1,\ldots,F_n) = F_1 \cup \ldots \cup F_n$ for each $(F_1,\ldots,F_n) \in \prod\{\mathcal{F}(\mathcal{V}_i) : i = 1,\ldots,n\}$. It is easy to verify that $u$ is continuous. Note that $u(\prod\{\mathcal{F}(\mathcal{V}_i) : i = 1,\ldots,n\}) = \left\langle \mathcal{V}\right\rangle_{\mathcal{F}(X)}$. Thus $\left\langle \mathcal{V}\right\rangle_{\mathcal{F}(X)}$ is connected. Now choose $V_S \in \mathcal{V}$ such that $x_S \in V_S$. According to Theorem \ref{SCC} the set $\mathcal{S}_c(V_S)$ is connected. It follows that $\mathcal{D} = \{T \cup F : T \in \mathcal{S}_c(V_S) \textnormal{ and } F \in \left\langle \mathcal{V}\right\rangle_{\mathcal{F}(X)}\}$, being a continuous image of $\mathcal{S}_c(V_S) \times (\left\langle \mathcal{V}\right\rangle_{\mathcal{F}(X)})$, is connected. On the other hand, it is easy to verify that $\mathcal{D}$ is dense in $\left\langle \mathcal{V}\right\rangle$. Therefore $\left\langle \mathcal{V}\right\rangle$ is a connected neighborhood of $S$ contained in $O$.
\endproof

\begin{theorem}
If $\mathcal{S}_c(X)$ is locally connected, then $X$ is locally connected.
\end{theorem}

\proof
Fix $x \in X$. Assume that $X$ is not locally connected at $x$. Fix an open neighborhood $U$ of $x$ in $X$ such that for any open set $V$ in $X$ with $x \in V \subseteq U$ the set $V$ is not connected. Note that $x$ is not an isolated point of $X$ and pick $S \in \mathcal{S}_c(X)$ such that $x = x_S$ and $S \subseteq U$.  Since $\mathcal{S}_c(X)$ is locally connected, we can find a connected open neighborhood $O$ of $S$ such that $S \in O \subseteq \left\langle U\right\rangle$. Let $N = \bigcup O$. It follows from Lemma \ref{LUA} that $N$ is open in $X$. Since $x \in N \subseteq U$, the set $N$ is not connected. Choose disjoint nonempty open subsets $V_0$ and $W_0$ of $X$ such that $N = V_0 \cup W_0$. For every $n \in \N^+$ construct disjoint nonempty open sets $V_n$ and $W_n$ with $x \in V_n \subseteq U$ recursively as follows. If $V_n$ has been constructed, since $x \in V_n \subseteq U$ we can find two disjoint nonempty open sets $V_{n+1}$ and $W_{n+1}$ such that $V_n = V_{n+1} \cup W_{n+1}$ and $x \in V_{n+1}$. Note that $N \setminus W_k = V_k \cup \bigcup\{W_n : n \in k\}$ is open, for each $k \in \N$.\medskip

\textbf{Claim 1.} The set $W_n$ is discrete for each $n \in \N$.\medskip

\textit{Proof of the Claim.} Assume on the contrary that $W_m$ is discrete for some $k \in \N$. Let $y$ be a non isolated point in $W_m$. Choose a sequence $S_0 \in O$ containing $y$ and a canonical open set $\left\langle \mathcal{V}\right\rangle$ of $\mathcal{S}_c(X)$ such that $S_0 \in \left\langle \mathcal{V}\right\rangle \subseteq O$. For each $V \in \mathcal{V}$ choose $y_V \in V$. Since $y$ is not isolated we may pick $S_y \in \mathcal{S}_c(X)$ such that $y \in S_y \subseteq \bigcup \mathcal{V}$ and $S_y$ converges to $y$. Let $S_1 = S_y \cup \{x_V : V \in \mathcal{V}\}$. Then $y \in S_1 \in \left\langle \mathcal{V}\right\rangle \subseteq O$ and $S_1$ converges to $y$. Since the family $\{W_n : n \in \N\}$ is disjoint and $y \in W_m$, we can fix $k \in \N$ such that $S_1 \cap W_k = \emptyset$. Pick $z \in W_k$. Since $W_k \subseteq N$ we may choose $T \in O$ such that $z \in T$. Then $T \cap W_k \not= \emptyset$. Note that the sets $O_1 = \left\langle N \setminus W_k\right\rangle$ and $O_2 = \left\langle W_k,N\right\rangle$ are open in $\mathcal{S}_c(X)$. Moreover, $O \subseteq O_1 \cup O_2$, $O_1 \cap O_2 = \emptyset$, $S_1 \in O \cap O_1$ and $T \in O \cap O_2$. Thus $O$ is disconnected, a contradiction. In this way, we have proved the Claim.\medskip

\textbf{Claim 2.} If $y \in N$ and $y$ is isolated in $X$, then $y \in T$ for each $T \in O$.\medskip

\textit{Proof of the Claim.} Fix $T \in O$. Suppose that $y \not \in T$. Since $y \in N$, we may choose a sequence $S_y \in O$ containing $y$. Consider the open sets $O_1 = \left\langle \{y\},N\right\rangle$ and $O_2 = \left\langle N \setminus \{y\}\right\rangle$. Note that $O \subseteq O_1 \cup O_2$, $O_1 \cap O_2 = \emptyset$, $S_y \in O \cap O_1$ and $T \in O \cap O_2$. It follows that $O$ is disconnected, a contradiction. Thus, $y \in T$ for each $T \in O$.\medskip

Let $W = \bigcup\{W_n : n \in \N\}$. Note that $W$ is an infinite subset of $N$ and $x \not \in W$. By Claim 1 the open set $W_n$ is discrete for each $n \in \N$, and as a consequence, each point of $W$ is isolated in $X$. It follows from Claim 2 that $W \subseteq T$ for each $T \in O$. In particular, $W \subseteq S$. Choose a canonical open set $\left\langle \mathcal{V}\right\rangle$ of $\mathcal{S}_c(X)$ such that $S \in \left\langle \mathcal{V}\right\rangle \subseteq O$. Since $\mathcal{V}$ is finite and $W$ is infinite, we can chose $V_0 \in \mathcal{V}$ such that $W \cap V_0$ is infinite. Pick $y \in W \cap V_0$ and let $T = S \setminus\{y\}$. Then $T \in \left\langle \mathcal{V}\right\rangle \subseteq O$, However $T$ does not contain $W$, contradicting the Claim 2. Thus $X$ is locally connected at the point $x$. Since $x$ is arbitrary, the space $X$ is locally connected.
\endproof

\begin{corollary}\label{SCLC}
The space $\mathcal{S}_c(X)$ is locally connected iff $X$ is locally connected.
\end{corollary}

It would be interesting to explore other connectedness like properties in the spirit of Corollaries \ref{SCCc} and \ref{SCLC}.

\section{path-wise connectedness on $\mathcal{S}_c(X)$}

In what follows, we shall prove that the path-wise connectedness of $\mathcal{S}_c(X)$ implies the pathwise connectedness of $X$.
We also deal with the number of path-wise connected components of $\mathcal{S}_c(X)$.
The following Lemma was proved in \cite{may-pat-rob} and will be a very important tool to establish the main result of the section.

\begin{lemma}\label{LMPR}
Let $X$ be a space. If $A$ is a compact connected subset of $\mathcal{K}(X)$ and $S \in A$, then each component of $\bigcup A$ intersects $S$.
\end{lemma}

The next theorem improves Theorem 4.11 from \cite{may-pat-rob}.

\medskip

Let $I$ denote the unit interval  $[0,1]$ of the real numbers $\R$.

\begin{theorem}\label{SCCT}
 If $\alpha : I \to \mathcal{S}_c(X)$ and $p \in \alpha(0)$, then there exists a continuous map $f : I \to X$ such that $f(0) = p$ and $f(t) \in \alpha(t)$ for all $t \in I$.
\end{theorem}

\proof
Let $K = \bigcup \alpha(I) \subseteq X$. The set $K$ is compact because of \cite[Theorem 2.5]{emi}. For each $n \in \N$ and  $i \in 2^n$, we  let $I_{n,i} = [i/2^n,(i+1)/2^n]$. Fix $n \in \N$.  For every $i \in 2^n$, define $p_{n,i}$ and $K_{n,i}$ recursively as follows. Let $p_{n,0} = p$ and $K_{n,0}$ be the connected component of $\bigcup\alpha(I_{n,0})$ which contains $p_{n,0}$. For $0 < i < 2^n$ , by Lemma \ref{LMPR}, we may choose $p_{n,i} \in K_{n,i-1} \cap \alpha(i/2^n)$. Let $K_{n,i}$ be a connected component of $\bigcup\alpha(I_{n,i})$ which contains $p_{n,i}$. We again apply  \cite[Theorem 2.5]{emi} to see that $K_{n,i}$ is compact for each $i \in 2^n$. Now, for each $n \in \N$ we consider the  compact and connected subspace $G_n = \bigcup\{I_{n,i} \times K_{n,i} : i \in 2^n\}$ of $I \times K$. Since  $\mathcal{CL}(I \times K)$ is compact, there exists a complete accumulation point $G$ of $\{G_n : n \in \N\}$ in $\mathcal{CL}(I \times K)$. By using the next three claims,  we will prove that $G$ is the graph of the promised map $f$.\medskip

\textbf{Claim 1.} $G \cap (\{t\} \times K) \subseteq \{t\} \times \alpha(t)$ for each $t \in I$.\medskip

\textit{Proof of the Claim.} By contradiction, assume that there exists $t \in I$ and a point $x \in K \setminus \alpha(t)$ such that $(t,x) \in G$. Choose disjoint open neighborhoods $U_x$ and $U_t$ of $x$ and $\alpha(t)$, respectively. By the continuity of $\alpha$, it is possible to find an open neighborhood $V$ of $t$ in $I$ and $N \in \N$ such that if $n \geq N$, $i \in 2^n$ and $V \cap I_{n,i} \not= \emptyset$. Then  we have that $\alpha(I_{n,i}) \subseteq \left\langle U_t \right\rangle$. Notice that $O  = \left\langle V \times U_x , I \times K\right\rangle$ is an open neighborhood of $G$ in $\mathcal{CL}(I \times K)$. Fix $n \geq N$. We will verify that $G_n \not \in O$. Choose $i \in 2^n$ and consider the following two cases: If $V \cap I_{n,i} = \emptyset$, then  it is clear that $(I_{n,i} \times K_{n,i}) \cap (V \times U_x) = \emptyset$. If $V \cap I_{n,i} \not= \emptyset$, then $\big(\bigcup\alpha(I_{n,i})\big) \cap U_x \subseteq U_t \cap U_x = \emptyset$. In particular, $K_{n,i} \cap U_x = \emptyset$ and as a consequence $(I_{n,i} \times K_{n,i}) \cap (V \times U_x) = \emptyset$. We have proved that, for each $n \in \N$,  the intersection  $G_n \cap (V \times U_x)$ is empty. Thus, we have that  $O$ is  a neighborhood of $G$ which intersect at most the first $N$ elements of $\{G_n : n \in \N\}$, which is impossible. Therefore, $G \cap (\{t\} \times K) \subseteq \alpha(t)$.\medskip

\textbf{Claim 2.} $\pi_I(G) = I$\footnote{ Here $\pi_I : I \times K \to I$ is the projection onto the first factor.}  and $(0,p) \in G$.

\medskip

\textit{Proof of the Claim.} First we will prove that $\pi_I(G) = I$. Assume on the contrary that $\pi_I(G) \not= I$ and choose a proper open subset $U$ of $I$ such that $\pi_I(G) \subseteq U$. Since $G \in \left\langle U \times K \right\rangle$, we can find $n \in \N$ such that $G_n \in \left\langle U \times K \right\rangle$. Then, $\pi_I(G_n) \subseteq U$, but this contradicts the fact  $\pi_I(G_n)= I$. Thus, we must have that  $\pi_I(G) = I$.

Now we will verify that $(0,p) \in G$. Suppose that $(0,p) \not\in G$ and choose $U = (I \times K) \setminus \{(0,p)\}$. Observe that $G \in \left\langle U\right\rangle$. By construction, we know that $(0,p) \in I_{n,0} \times K_{n,0} \subseteq G_n$ for each $n \in \N$; that is, $G_n \not \in \left\langle U\right\rangle$. As a consequence $\left\langle U\right\rangle \cap \{G_n : n \in \N\} = \emptyset$, a contradiction. This shows that $(0,p) \in G$.\medskip

\textbf{Claim 3.} $G$ is the graph of a map.\medskip

\textit{Proof of the Claim.} Let $t \in I$ and assume that $(t,x_1), (t,x_2) \in G$. We must prove that $x_1 = x_2$. Assume that $x_1 \not= x_2$. In virtue of the first claim, we must have that $x_1,x_2 \in \alpha(t)$. Let $\mathcal{U} = \{U_1,U_2\}$ be a disjoint open cover of $\alpha(t)$ such that $x_1 \in U_1$ and $x_2 \in U_2$. By the continuity of $\alpha$, we can find an open interval $V$ in $I$ containing $t$ and $N \in \N$ such that if $n \geq N$, $i \in 2^n$ and $V \cap I_{n,i} \not= \emptyset$, then $\alpha(I_{n,i}) \subseteq \left\langle \mathcal{U} \right\rangle$. Note that $G \in O := \left\langle V \times U_1 ,V \times U_2, I \times K\right\rangle$. We will verify that $G_n \not \in O$ whenever $n \geq N$. Fix $n \geq N$. Now, consider the sets $G_n^\prime = \{ I_{n,i} \times K_{n,i} : i \in 2^n  \textnormal{ and } I_{n,i} \cap V = \emptyset\}$ and $G_n^{\prime\prime} = \{ I_{n,i} \times K_{n,i} : i \in 2^n  \textnormal{ and } I_{n,i} \cap V \not= \emptyset\}$. In one hand, it is clear that $G_n^\prime \cap (V \times U_i) = \emptyset$ for $i = 1,2$. On the other hand, $G_n^{\prime\prime}$ is connected by the construction. Moreover, by the election of $V$, $G_n^{\prime\prime} \subseteq (I \times U_1) \cup (I \times U_2)$. So there exists $i_0 \in \{1,2\}$ such that $G_n^{\prime\prime} \cap (V \times U_{i_0}) = \emptyset$. It follows that $G_n \cap (V \times U_{i_0}) = \big(G_n^{\prime} \cap (V \times U_{i_0})\big) \cup \big(G_n^{\prime\prime} \cap (V \times U_{i_0})\big) = \emptyset$. This shows that $G_n \not \in O$. It follows that $O$ is  a neighborhood of $G$ which intersect at most the first $N$ elements of $\{G_n : n \in \N\}$, a contradiction. Thus $x_1 = x_2$.

\medskip

Our required function is $f = G$. It follows from Claims 1, 2 and 3 that $f$ is a function from $[0,1]$ into $K$ such that $f(0) = p$ and $f(t) \in \alpha(t)$ for each $t \in I$. We know that the graph $G$ of $f$ is closed and $K$ is compact, so we may apply the Closed Graph Theorem (see \cite[Ex. 3.1.D (a)]{eng}) to see that $f$ is continuous.
\endproof

\begin{corollary}
If $\mathcal{S}_c(X)$ is path-wise connected, then $X$ is path-wise connected.
\end{corollary}

\proof
Fix two distinct points $p, q \in X$.  Since $X$ is connected, see  Corollary \ref{SCCc}, the space $X$ cannot have isolated points. Choose disjoint sequences $S_p,S_q \in \mathcal{S}_c(X)$ converging to $p$ and $q$, respectively. By Theorem \ref{SCCT}, there exists a path $f_1 : [0,1/3] \to \mathcal{S}_c(X)$ such that $f_1(0) = p$ and $r := f_1(1/3) \in S_q$. Choose a nontrivial sequence $S_r$ in $\bigcup f_1([0,1/3))$ converging to $r$. We again apply Theorem  \ref{SCCT} to find a path $f_3 : [2/3,1] \to \mathcal{S}_c(X)$ such that $f_3(1) = q$ and $s := f_3(2/3) \in S_r$. Now note that there exists a path $f_2 : [1/3,2/3] \to \mathcal{S}_c(X)$ so that $f_2(1/3) = r$ and $f_2(2/3) = s$. Then, $f = f_1 \cup f_2 \cup f_3$ is a path that connects  $p$  with $q$.
\endproof

The above result is a very interesting relation between  $\mathcal{S}_c(X)$ and $X$. However,  we do not have  any positive result in the opposite direction.

\begin{problem}
Give conditions on $X$ under which $\mathcal{S}_c(X)$ must be path-wise connected.
\end{problem}

It is known that there are path-wise connected spaces $X$ for which the space $\mathcal{S}_c(X)$ is not connected \cite[Ex. 2.8]{sal-yas}. However, we do not know so much about the number of path-wise connected components  of $\mathcal{S}_c(X)$. Next we shall give two examples of spaces $X$ for which their hyperspace $\mathcal{S}_c(X)$ has $\mathfrak{c}$-many path-wise connected components.

\medskip

We need the next lemma which was essentially proved in \cite[Lemma 2.7]{sal-yas}.

\begin{lemma}\label{LYS}
Let $f : [0,1] \to \mathcal{S}_c(X)$ be a continuous function. Suppose $C$ is a nonempty closed subset of $\bigcup f([0, 1])$. If $s = \sup \{t \in [0, 1] : f(t) \cap C \not= \emptyset\}$, then $f(s) \cap C \not= \emptyset$.
\end{lemma}

The first example is a dendroid that also was considered in \cite[Example 4.7]{may-pat-rob}.

\begin{example}
There exists a dendroid $X \subset \R^2$ such that $\mathcal{S}_c(X)$ has $\mathfrak{c}$-many path-wise connected components.
\end{example}

\proof
 For $x, y \in \R^2$ we let $[x,y] = \{ty + (1-t)x : t \in [0,1]\}$ and $(x,y] = \{ty + (1-t)x : t \in (0,1]\}$. Set $u = (0,1)$, $v = (0,0)$, $w = (0,-1)$ and $x_n = (1/n,0)$ for each $n \in \mathbb{Z} \setminus \{0\}$. Consider the dendroid
$$X = [u,w] \cup \bigcup_{n \in {\mathbb{Z}^+}}[u,x_n] \cup \bigcup_{n \in {\mathbb{Z}^-}}[w,x_n]$$
endowed with the topology inherited of $\R^2$. We remark that  the dendroid $\mathcal{S}_c(X)$ has at most $\mathfrak{c}$-many path-wise connected components. For any infinite set $A \subseteq {\mathbb{Z}^+}$ let us define $S_A = \{v\} \cup\{x_n : n \in A \cup {\mathbb{Z}^-}\} \in \mathcal{S}_c(X)$.\medskip

\textbf{Claim.} If $A$ and $B$ are two almost disjoint subsets of ${\mathbb{Z}^+}$, then there is not a path between $S_A$ and $S_B$ in $\mathcal{S}_c(X)$.\medskip

\textit{Proof of the Claim.} Assume on the contrary that there is a path $f : [0,1] \to \mathcal{S}_c(X)$ such that $f(0) = S_A$ and $f(1) = S_B$. Let $X_n = \bigcup\{(u,x_k] : k \in B \textnormal{ and } k \geq n\}$ if $n \in B$ and let $X_n = \bigcup\{(w,x_k] : k \in {\mathbb{Z}^-} \textnormal{ and } k \leq n\}$ if $n \in {\mathbb{Z}^-}$. Consider the nonempty set $C = \{t \in [0,1] : \exists n \in B \cup {\mathbb{Z}^-}(f(t) \cap X_n = \emptyset)\}$ and  $s := \sup C$. Let $\mathcal{U}$ be a disjoint family of open subsets of $X$ with diameter less than $1/2$, such that $f(s) \in \left\langle \mathcal{U}\right\rangle$. In order finish the proof of the Claim, we consider two cases (in both cases we will get a contradiction):

\smallskip

\textit{Case 1.} $s \not\in C$. Note that $0 < s$. Choose $r \in (0,s) \cap C$ such that $[r,s] \in f^{-1}(\left\langle \mathcal{U}\right\rangle)$. Then we can find $n_0 \in B \cup {\mathbb{Z}^-}$ for which  $f(r) \cap X_{n_0} = \emptyset$. Without lost of generality, we can suppose that $n_0 \in B$.  Since $s \not\in C$ it then follows  that $f(s) \cap X_n \not= \emptyset$ for each $n \in B \cup {\mathbb{Z}^-}$. As a consequence $x_{f(s)} = v$. Choose $U_0 \in \mathcal{U}$ such that $v \in U_0$. Since $v \in U_0$ and $f(s) \cap X_n \not= \emptyset$ for each $n \in B$, we can find $k \in B$ such that $k \geq n_0$ and $f(s) \cap U_0 \cap (u,x_k] \not= \emptyset$. Applying Theorem 2.5 from \cite{emi}, we obtain that the set $K = \bigcup f([r,s])$ is compact. Since $K \subseteq \bigcup \mathcal{U}$, by the election of $\mathcal{U}$, the set $K_0 = K \cap U_0 \cap (u,x_k] = K \cap \textnormal{cl}(U_0) \cap [u,x_k]$ is compact. Besides, since $(u,x_k]$ is open in $X$, the set $V = (u,x_k] \cap U_0$ is an open neighborhood of $K_0$ in $X$. Let $\mathcal{V} = \mathcal{U} \cup \{V\}$. Consider the nonempty set $D = \{q \in [r,s] : f(q) \cap K_0 \not= \emptyset\}$ and  $i := \inf D$. We apply the infimum version of Lemma \ref{LYS} to see that $i \in D$. Then,  we have that $r < i \leq s$ and $f(i) \in \left\langle \mathcal{V} \right\rangle$. However, for each $r < t < i$ observe that $f(t) \not \in \left\langle \mathcal{V} \right\rangle$, contradicting the continuity of $f$.

\smallskip

\textit{Case 2.} $s \in C$. In this case we must have that $s < 1$. Choose $n_0 \in B \cup {\mathbb{Z}^-}$ such that $f(s) \cap X_{n_0} = \emptyset$.  We may suppose, without lost of generality, that $n_0 \in {\mathbb{Z}^-}$. Choose $t \in (s,1)$ such that $[s,t] \in f^{-1}(\left\langle \mathcal{U}\right\rangle)$. The condition $t \not\in C$ implies that $f(t) \cap X_n \not= \emptyset$ for each $n \in B \cup {\mathbb{Z}^-}$. As a consequence, we obtain that $x_{f(t)} = v$. Choose $U_0 \in \mathcal{U}$ such that $v \in U_0$. Since $v \in U_0$ and $f(t) \cap X_n \not= \emptyset$ for each $n \in {\mathbb{Z}^-}$,  there exits $k \in {\mathbb{Z}^-}$ such that $k \leq n_0$ and $f(t) \cap U_0 \cap (u,x_k] \not= \emptyset$. As above,  the sets $K = \bigcup f([s,t])$ and $K_0 = K \cap U_0 \cap (u,x_k] = K \cap \textnormal{cl}(U_0) \cap [u,x_k]$ are compact. Moreover, the set $V = (u,x_k] \cap U_0$ is an open neighborhood of $K_0$ in $X$. Consider the family $\mathcal{V} = \mathcal{U} \cup \{V\}$. Let $D = \{q \in [s,t] : f(q) \cap K_0 \not= \emptyset\}$ and $i := \inf D$. By the infimum version of Lemma \ref{LYS}, $i \in D$. Note that $s < i \leq t$ and $f(i) \in \left\langle \mathcal{V} \right\rangle$. Finally, for $s < r < i$ we have that $f(r) \not \in \left\langle \mathcal{V} \right\rangle$, which contradicts the continuity of $f$.\medskip

Let $\mathcal{A}$ be an almost disjoint family of subsets ${\mathbb{Z}^+}$ of size $\mathfrak{c}$ (for the existence of such a family see 6Q in \cite{gaj}). It follows from the Claim that for  distinct elements $A,B \in \mathcal{A}$ the sequences $S_A$ and $S_B$ belong to distinct path-wise connected components. Therefore,   $\mathcal{S}_c(X)$ has $\mathfrak{c}$-many path-wise connected components.
\endproof

\begin{example}
The space of non-trivial convergent sequences on the Warsaw circle has $\mathfrak{c}$-many path-wise connected components.
\end{example}

\proof
The Warsaw circle is the subspace $X = X_1 \cup X_2 \cup X_3$ of $\R^2$, where $X_1 = \{(a,\left|\sin(\pi/a)\right|) : a \in (0,1]\}$, $X_2 = \{(a,b) \in \R^2 : (a-1/2)^2 + b^2 = 1/4 \textnormal{ and } b < 0\}$ and $X_3 = \{(0,b) : b \in [0,1]\}$. Consider the subspace $D = X \cap (\R \times \{0\})$. For each two points $x,y \in X$ define the natural number $\rho(x,y) = \left|I \cap D\right|$ where $I$ is the only arc in $X$ with endpoints $x$ and $y$. For $x \in X$ and $S \in \mathcal{S}_c(X)$ we let $\rho(x,S) = \min\{\rho(x,y) : y \in S\}$. Now for $S_1,S_2 \in \mathcal{S}_c(X)$ we define
$$
\rho(S_1,S_2) = \max\{\sup\{\rho(x,S_2) : x \in S_1\},\sup\{\rho(y,S_1) : y \in S_2\}\}.
$$
 Let $Y  = \{S \in \mathcal{S}_c(X) : \left|S \cap X_1\right| = \omega \textnormal{ and } x_S \in X_3\}$. Given $S \in Y$ we say that a canonical neighborhood $\left\langle \mathcal{U}\right\rangle$ of $S$ in $\mathcal{S}_c(X)$ is a \textit{suitable} neighborhood of $S$ if for every $U \in \mathcal{U}$ the closure of each component of $U \cap X_1$ is an arc and contains at most one point of $D$.
\medskip

\textbf{Claim 1.} If $f : [0,1] \to Y$ is a path,  $\left\langle \mathcal{U}\right\rangle$ is a suitable neighborhood of $f(1)$ and $f([0,1]) \subseteq \left\langle \mathcal{U}\right\rangle$, then $\rho(f(0),f(1)) < \omega$.\medskip

\textit{Proof of the Claim.} Assume on the contrary that $\rho(f(0),f(1)) = \omega$. Consider the compact space  $K = \bigcup f([0,1])$.  In each one of the following cases we are gong to get a contradiction, so the Claim 1 will be proved.

\textit{Case 1.} $\sup\{\rho(x,f(1)) : x \in f(0)\} = \omega$. Fix $x_0 \in f(0) \cap X_1$ such that $\rho(x_0,f(1)) \geq 2$. Choose $U \in \mathcal{U}$ such that $x_0 \in U$, let $U_0$ be the connected component of $U$ which contains $x_0$ and let $K_0 = K \cap U_0$. Note that since the elements of $\mathcal{U}$ are pairwise  disjoint and $K \subseteq \bigcup \mathcal{U}$, the set $K_0$ is compact. Let $A = \{a \in [0,1] : f(a) \cap K_0 \not = \emptyset\}$. By Lemma \ref{LYS}, $s := \sup A \in A$. Since $x_0 \in U_0$ and $\textnormal{cl}(U_0)$ is an arc with at most one point of $D$, it follows from $\rho(x_0,f(1)) \geq 2$ that we must have $f(1) \cap U_0 = \emptyset$. Hence, $s < 1$. Let $\mathcal{V} = \mathcal{U} \cup \{U_0\}$. Then $f(s) \in \left\langle \mathcal{V}\right\rangle$. By continuity, we may choose $t \in (s,1)$ such that $f([s,t]) \subseteq \left\langle \mathcal{V}\right\rangle$. Then $f(t) \cap U_0 \not=\emptyset$ and so $f(t) \cap K_0 = f(t) \cap U_0 \not= \emptyset$. It follows that $t \in A$, a contradiction.

\smallskip

\textit{Case 2.} $\sup\{\rho(y,f(0)) : y \in f(1)\} = \omega$. Pick $y_0 \in f(1) \cap X_1$ for which $\rho(y_0,f(0)) \geq 2$. Consider $U \in \mathcal{U}$ such that $y_0 \in U$ and denote  by $U_0$ the connected component of $U$ which contains $y_0$. As above $K_0 = K \cap U_0$ is compact. Let us consider the set $A = \{a \in [0,1] : f(a) \cap K_0 \not = \emptyset\}$. By the infimum version of Lemma \ref{LYS}, $i := \inf A \in A$. We know that $y_0 \in U_0$, $\textnormal{cl}(U_0)$ is an arc with at most one point of $D$,  and $\rho(y_0,f(0)) \geq 2$ which imply $f(0) \cap U_0 = \emptyset$. It follows that $i > 0$. Consider the family $\mathcal{V} = \mathcal{U} \cup \{U_0\}$ and note that $f(i) \in \left\langle \mathcal{V}\right\rangle$. Select $j \in (0,i)$ such that $f([j,i]) \subseteq \left\langle \mathcal{V}\right\rangle$. Note that $f(j) \cap U_0 \not=\emptyset$ and so $f(j) \cap K_0 = f(j) \cap U_0 \not= \emptyset$. Therefore, $j \in A$, a contradiction.
\medskip

\textbf{Claim 2.} If $S_0,S_1 \in Y$ and $\rho(S_0,S_1) = \omega$, then there is no a path between $S_0$ and $S_1$ in $\mathcal{S}_c(X)$.\medskip

\textit{Proof of the Claim.} Suppose that $f : [0,1] \to \mathcal{S}_c(X)$ is a path with $f(0) = S_0$ and $f(1) = S_1$. A similar argument to the one  given in the proof of \cite[Example 2.8]{sal-yas} shows that there is no a path between a sequence in $Y$ and a sequence in $\mathcal{S}_c(X) \setminus Y$. So, we must have that $f([0,1]) \subseteq Y$. Let $A = \{a \in [0,1] : \rho(f(0),f(a)) < \omega\}$ and  $s = \sup A$. Let $\left\langle \mathcal{U}\right\rangle$ be a suitable neighborhood of $f(s)$. We will prove that $s \not\in A$. Suppose that $s \in A$. Then we must have $s < 1$. Choose $t \in (s,1)$ such that $f([s,t]) \subseteq \left\langle \mathcal{U}\right\rangle$. Claim 1 implies that $\rho(f(s),f(t)) < \omega$. By our assumption, $\rho(f(0),f(s)) < \omega$. So we must have $\rho(f(0),f(t)) < \omega$ and hence $t \in A$, but this is a contradiction. This shows that $s \not\in A$. Then $s > 0$. Choose $r \in (0,s)$ such that $r \in A$ and $f([r,s]) \subseteq \left\langle \mathcal{U}\right\rangle$. Claim 1 implies that $\rho(f(r),f(s)) < \omega$. Since $\rho(f(0),f(r)) < \omega$, we obtain that $\rho(f(0),f(s)) < \omega$. Hence $s \in A$, a contradiction.

\medskip

To finish the proof let $\mathcal{A}$ be an almost disjoint family of subsets $\N^+$ of cardinality $\mathfrak{c}$. Let $k_n = \sum_{k = 1}^n k$ and $a_n = 1/k_n$ for each $n \in \N^+$.  Note that $\rho((a_n,0),(a_{n+1},0)) \geq n$ for all $n \in \N^+$. For each $A \in \mathcal{A}$ let $S_A = \{(0,0),(a_n,0) : n \in A\}$. Choose $A,B \in \mathcal{A}$ such that $A \not= B$. By the construction, we have that $\rho(S_A,S_B) = \omega$. It follows from Claim 2 that $S_A$ and $S_B$ belong to distinct path-wise connected components of $\mathcal{S}_c(X)$. Therefore, $\mathcal{S}_c(X)$ admits  $\mathfrak{c}$-many path-wise connected components.
\endproof

We do not know if the behavior of the number of path-wise connected components is always as in the above examples.

\begin{question}
Let $X$ be a path-wise connected space. Can $\mathcal{S}_c(X)$ have a finite or even a countable number of path-wise connected components?
\end{question}



\begin{thebibliography}{99}
\bibitem{arh} A. V.  Arhangel'skii and A. M.  Al Shumrani,
\newblock {\it Jointly partially metrizable spaces},
\newblock {C. R. Acad. Bulgare Sci.}  \textbf{65}  (2012),  no. 6, 727--732.

\bibitem{eng} R. Engelking,
\newblock {\it  General topology},
\newblock {Sigma Series in Pure Mathematics, vol. 6, Heldermann Verlag, Berlin, 1989, Translated from the Polish by the author.}

\bibitem{sal-yas} S. Garcia-Ferreira and Y. F. Ortiz-Castillo,
\newblock {\it The hyperspace of convergent sequences},
\newblock to appear in {Top. Appl.}.

\bibitem{gaj} L. Gillman and M. Jerison,
\newblock {\it Rings of Continuous Functions},
\newblock  {Grad. Texts in Math. 43, Springer Verlag, 1976.}

\bibitem{gu} V. Gutev,
\newblock {\it Selections and hyperspaces},
\newblock  {Recent progress in general topology. III,  535--579, Atlantis Press, Paris, 2014.}

\bibitem{may-pat-rob} D. Maya, P. Pellicer-Covarrubias and R. Pichardo-Mendoza,
\newblock {\it On the hyperspace of convergent sequences},
\newblock submitted.

\bibitem{emi} E. Michael,
\newblock {\it Topologies on spaces of subsets},
\newblock {Trans. Amer. Math. Soc., \textbf{71} (1951), 152--182.}

\bibitem{jvpp} J. van Mill, J. Pelant and R. Pol,
\newblock {\it Selections that characterize topological completeness},
\newblock {Fund. Math.},  \textbf{149}  (1996),  no. 2, 127–-141.

\bibitem{peng} P. Liang-Xue and  G. Zhi-Fang,
\newblock {\it A note on joint metrizability of spaces on families of subspaces},
\newblock {Top. Appl.},  \textbf{188}  (2015), 1–-15.
\end{thebibliography}
\end{document}